\documentclass[journal]{IEEEtran}

\ifCLASSINFOpdf
 \else
 \fi

\usepackage{multirow}
\usepackage{cite}
\usepackage{amsmath,amssymb,amsfonts}
\usepackage[dvipsnames]{xcolor}
\usepackage{algorithm,algorithmic} 
\usepackage{amssymb}
\usepackage{mathrsfs}
\usepackage{float,graphicx}
\usepackage{mathtools}
 \usepackage{lscape}
\def\mathclap#1{\text{\hbox to 0pt{\hss$\mathsurround=0pt#1$\hss}}}
\usepackage[mathlines]{lineno} 
\newcommand{\qedsymbol}{\rule{2mm}{2mm}}
\begin{document}  

\title{Optimal Coordination of Platoons of Connected and Automated Vehicles at Signal-Free Intersections} 

\author{Sharmila~Devi~Kumaravel,~\IEEEmembership{Student Member,~IEEE,} Andreas~A.~Malikopoulos,~\IEEEmembership{Senior Member,~IEEE}
        Ramakalyan~Ayyagari,~\IEEEmembership{Senior Member,~IEEE}
       % and~Jane~Doe,~\IEEEmembership{Life~Fellow,~IEEE}% <-this % stops a space 
\thanks{This research was supported in part by National Institute of Technology, Tiruchirappalli, India, and in part by ARPAE's NEXTCAR program under the award number DE-AR0000796 and by the Delaware Energy Institute (DEI).}
\thanks{S. D. Kumaravel and R. Ayyagari are with Department of Instrumentation and Control Engineering, National Institute of Technology, Tiruchirappalli, India (e-mail: \texttt{info2sj@gmail.com}; \texttt{rkalyn@nitt.edu} ).}% <-this % stops a space
\thanks{A. A. Malikopoulos is with the Department of Mechanical Engineering, University of Delaware, Newark, DE 19716 USA (e-mail: \texttt{andreas@udel.edu}).}
\thanks{ }% <-this % stops a space
\thanks{ }}%August 30, 2018; revised october 26, 2018}}

\maketitle 

\begin{abstract}
	
In this paper, we address the problem of coordinating platoons of connected and automated vehicles crossing a signal-free intersection.  We present a decentralized, two-level optimal framework to coordinate the platoons with the objective to minimize travel delay and fuel consumption of every platoon crossing the intersection. At the upper-level, each platoon leader derives a proven optimal schedule to enter the intersection. At the low-level, the platoon leader derives their optimal control input (acceleration/deceleration) for the optimal schedule derived in the upper-level.  We validate the effectiveness of the proposed framework in simulation and show significant improvements both in travel delay and fuel consumption compared to the baseline scenarios where platoons enter the intersection based on  first-come-first-serve and longest queue first - maximum weight matching scheduling algorithms.   
\end{abstract}


\IEEEpeerreviewmaketitle
\section{Introduction}
\subsection{Motivation}
\IEEEPARstart{T}{raffic} congestion has become a severe issue in urban transportation networks across the globe. Transportation networks will account for nearly 70\% of travel in the world with more than 3 billion vehicles by 2050 \cite{iea2013policy}. The exponential growth in the number of vehicles and rapid urbanization have contributed to the steadily increasing problem of traffic congestion. The drivers lose 97 hours due to congestion and the cost of congestion was estimated to be \$87 billion a year, i.e., an average of \$1,348 per driver in US \cite{reed2019inrix}. Urban intersections in conjunction with the driver's response to various disturbances can aggravate congestion.  Efficient intersection control algorithms can improve mobility, safety and alleviate the severity of congestion and accidents. Recent advancements in vehicle-to-infrastructure and vehicle-to-vehicle (V2V) communication provide promising opportunities for control algorithms to reduce delay, travel time, fuel consumption, and emissions of vehicles \cite{rios2016survey}. The advent of connected and automated vehicles (CAVs) along with communication technologies can  enhance urban mobility with better options to travel efficiently \cite{zhao2019enhanced}. Moreover,  real-time information from CAVs related to their position, speed and acceleration through on-board sensors and V2V communication makes it possible to develop effective control algorithms for coordinating CAVs aimed at improving mobility and alleviate congestion.

\subsection{Related Work}
Several research efforts have proposed centralized and decentralized control algorithms for coordinating CAVs at intersections. Dresner and Stone \cite{dresner2008multiagent} presented a reservation scheme as an alternative approach to traffic lights for coordinating CAVs at an intersection.  Following this effort,  several centralized approaches have been reported in the literature to coordinate CAVs at signal-free intersections and other traffic scenarios, e.g., merging roadways \cite{huang2012assessing, lee2012development,  wu2013cooperative,rios2016automated,qian2017autonomous,lin2017autonomous,Rios2018}.  
Recently, Hart \emph{et al}. \cite{hart2019fail} developed an intersection control algorithm that considers safety and parametric uncertainties. Other research efforts in the literature have proposed decentralized control algorithms for coordinating vehicles at signal-free intersections. Wu \emph{et al}.\cite{wu2014distributed} proposed a decentralized control algorithm based on the estimated arrival time of CAVs at an intersection  without  eliminating stop-and-go  driving. To eliminate stop-and-go driving and minimize energy consumption, a  decentralized optimal control framework was presented  to coordinate CAVs at an intersection and for a corridor with different transportation scenarios\cite{malikopoulos2018decentralized, Malikopoulos2019CDC, Malikopoulos2020, zhao2018decentralized}. 

\par The road capacity and operational efficiency of the intersection can be increased significantly if the vehicles cross the intersection as platoons instead of crossing one after the other \cite{lioris2017platoons}. Prior to this, various research efforts in the literature address vehicle platooning at highways to increase fuel efficiency, traffic flow, comfort of driver, and safety.  Bergenhem \emph{et al}. \cite{bergenhem2012overview} presented a  detailed discussion on various research efforts in vehicle platooning systems at highways.  Vehicle platooning is not only beneficial at highways but also at urban traffic intersections.
Several research efforts have presented control algorithms in the literature to coordinate platoons at intersections.  
Jin \emph{et al}.  \cite{jin2013platoon} presented an intersection management under a multiagent framework in which the platoon leaders send the arrival time of vehicles and request to cross the intersection based on first-come-first-serve (FCFS) policy.  A hierarchical intersection management system was presented in  \cite{tallapragada2015coordinated} with the objective to minimize cumulative travel time and energy usage.  The research effort in  \cite{vial2016scheduling} proposed polynomial time algorithms to find schedules for the intersection with two-way traffic.  Bashiri  and   Fleming \cite{bashiri2017platoon} presented a control algorithm  that performed an extensive search among $\displaystyle n!$ schedules (which shoot up exponentially) to find a schedule with minimum average delay for $\displaystyle n$ platoons. Later, a greedy algorithm was presented in \cite{bashiri2018paim} to find the best schedule from all possible schedules that minimizes the total delay. The research efforts in the literature developed rule based control algorithm \cite{du2018v2x}, nonlinear control
algorithm \cite{di2019design}, and polling based control algorithm \cite{miculescu2019polling} for CAVs to gain access into the intersection. 

Recently, Feng \emph{et al}. \cite{feng2019composite} presented a reinforcement learning based control algorithm to plan the trajectories of platoons with the objective to maximize the throughput of signalized intersections. In order to  completely acquire the benefits of vehicle platooning at the intersections, effective scheduling and planning of the platoons are very essential.  Scheduling theory can offer  effective solutions to schedule the platoons at intersections. Various techniques in scheduling theory efficiently allocate limited resources to several tasks to optimize the performance measures. Scheduling theory based control algorithm is presented in \cite{giridhar2006scheduling} to coordinate vehicles in the urban roads. Li \emph{et al}.  \cite{li2006cooperative} presented a  safe driving for vehicle pairs to avoid collisions at the intersections. Schedule-driven control algorithms have been reported in the literature to evacuate all vehicles in minimum time at an intersection \cite{wu2009intersection} and for multiple intersections \cite{yan2011scheduling}.  
 A least restrictive supervisor was designed in \cite{colombo2014least} and \cite{ahn2016semi} to determine set of control actions for the vehicles to safely cross the intersection. The research effort in \cite{chalaki2019optimal}  derived the optimal schedule for CAVs and presented a closed-form analytical solution to derive optimal control input for vehicles at intersections. 

 Our proposed approach aims to overcome the limitations of existing approaches in the literature in the following ways:
 \begin{enumerate}
     \item The majority of the papers in the literature have proposed approaches for coordinating CAVs to cross the intersection one after another rather than platoons. The communication burden is significantly reduced when an intersection manager communicates only with the platoon leader instead of communicating with every CAV. Moreover, the capacity of the intersection significantly increases by vehicle platooning than allowing them to pass one after another.
     \item Most research efforts have employed a centralized approach for coordinating platoons of CAVs at an intersection. The approach is centralized if there is at least one task in the system that is globally decided for all vehicles by a single central controller. The decision that includes all vehicles will typically result in high communication and computational load. Furthermore, centralized approaches are ineffectual in handling single point failures.  On the other hand, decentralized approach reduces the communication requirements and are computationally efficient.      
 \end{enumerate}
   
We present a decentralized control framework for coordinating platoons of CAVs where each platoon leader communicates with other platoon leaders and a coordinator to derive the optimal schedule to cross the intersection. Furthermore, each platoon leader derives its optimal control input to cross the intersection while minimizing travel delay and fuel consumption.  

\subsection{Contributions of the paper}
%platoon coordinaton % decentralized
 The main contributions of the paper are the following. We present a decentralized, two-level optimal control framework to coordinate the platoons at an intersection. In the upper-level, we propose a proven optimal framework where each platoon computes the optimal schedule to minimize the travel delay of platoons. In the low-level, we present a closed-form analytical solution that provides the optimal control input  to minimize fuel consumption of vehicles.
 
\subsection{Organization of the paper}
The paper is organized as follows. In Section II, we formulate the problem, introduce the modeling framework, and present the upper-level framework that provides the optimal schedule for platoons. In Section III, we provide a closed-form, analytical solution of the low-level optimal control problem. In Section IV, we validate the effectiveness of the proposed optimal framework using VISSIM-MATLAB environment and present the simulation results. We conclude and discuss the potential directions for future work in Section V. 
\section{Problem Formulation}
\begin{figure*}[htbp!]
	\centering
	\includegraphics[height=4in,width=5in]{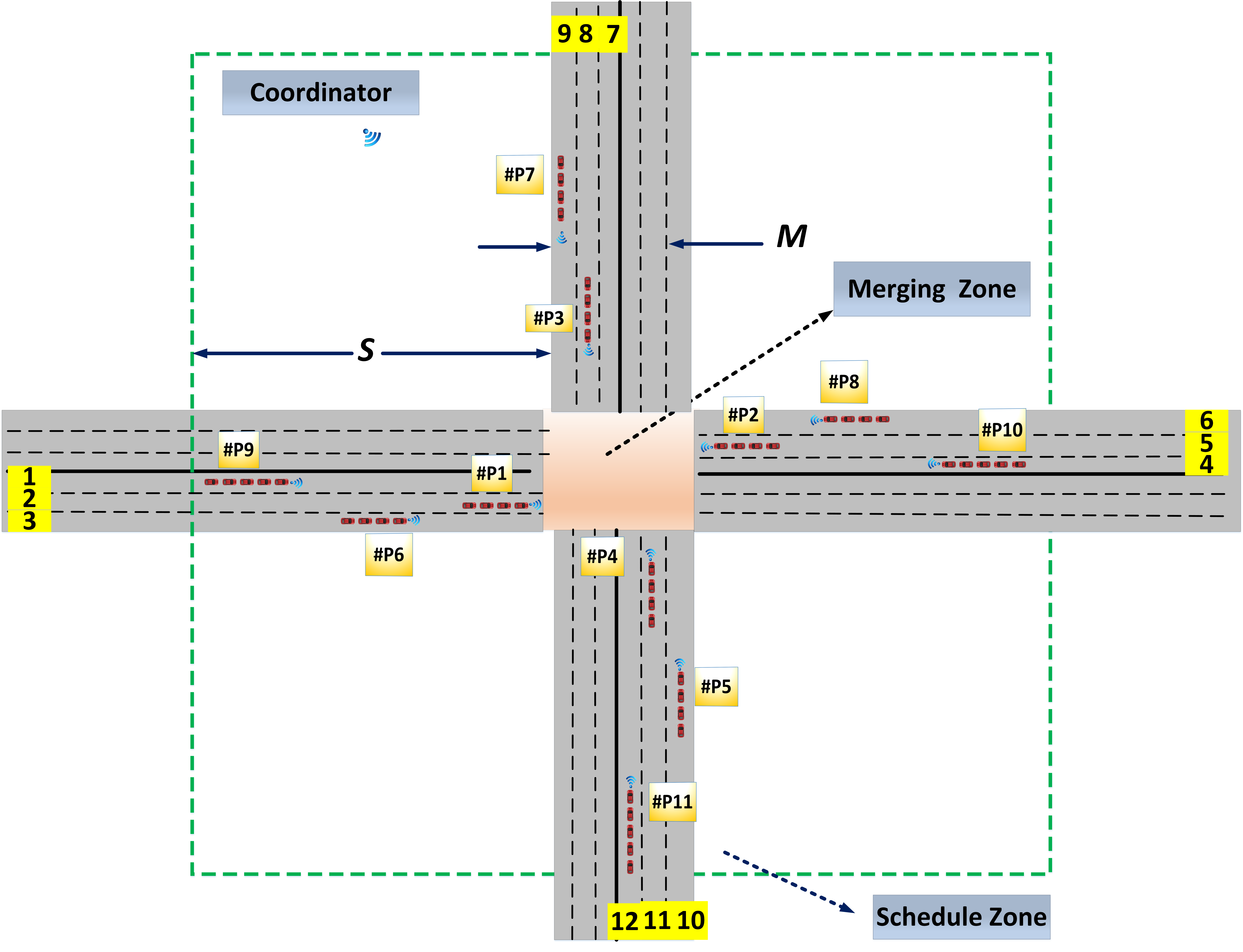}
	\caption{A signal-free intersection with platoons of CAVs.}
	\label{platoonInt}
\end{figure*}
We consider a signal-free traffic intersection (Fig. \ref{platoonInt}) for coordinating platoons of CAVs with minimum travel delay. The region at the center of the intersection is called \emph{merging zone}, which is the conflict area where potential lateral collisions of CAVs are possible. Although this is not restrictive, we consider the merging zone to be a square of side $\displaystyle M.$  The intersection has a \emph{schedule zone} and a coordinator that can communicate with the vehicles traveling inside the schedule zone. The distance from the entry point of the schedule zone until the entry point of the merging zone is $\displaystyle S$. The value of $\displaystyle S$ depends on the communication range capability of the coordinator.  The coordinator stores the information about the geometry and topology of the intersection. In addition, the coordinator stores the information about position, speed, acceleration/deceleration, and path of the platoons. Note that the coordinator acts as a database and does not take part in any of the decision making process. Each platoon leader can communicate with the coordinator, their followers and the other platoon leaders inside the schedule zone. 

 Let $\mathcal N(t)=\{1,\ldots N(t)\}$, $\displaystyle N(t) \in  \mathbb{N}$, be the queue of platoons inside the schedule zone. At the entry of the schedule zone, each platoon $j\in\mathcal N(t)$ leader broadcasts the information of the platoon to the coordinator and other platoon leaders. This information is the 6-tuple $\{n_j$, $N_{link}$, $N_{lane}$, $D_j$, $p_j$, $v_j\}$ in which $n_j$ denotes the number of vehicles in the platoon, $N_{link}$ denotes link number (link is the incoming road at the intersection), $N_{lane}$ denotes  lane number, $D_j$ denotes   routing decision, $p_j$ denotes  current position, and $v_j$ denotes current speed of the platoon. Based on the information from the coordinator and other platoon leaders, each platoon leader derives the time to enter the merging zone and optimal control input to cross the intersection. Each platoon leader broadcasts the schedule and optimal control input to the followers in the platoon, and then communicates the schedule to the coordinator.  The coordinator broadcasts the schedule of platoons inside the schedule zone to the leaders of the platoons entering the schedule zone.  \\\\
 In our modeling framework, we impose the following assumptions:\\\\
 \textbf{Assumption 1:} 
 There is no delay and communication errors between platoon leaders, the followers, and the coordinator.\\\\
 The assumption may be strong, but it is relatively straightforward to relax it as long as the measurement noise and delays are bounded \cite{lucas2020scalability} in a statistical sense.\\\\
\textbf{Assumption 2:}
The CAVs within the communication range form stable platoons, i.e., all the vehicles in the platoon
move at a consensual speed and maintain the desired space
between vehicles \cite{levine1966optimal}. \\ \\
Our primary focus is to coordinate the platoons of CAVs rather than the formation and stability of platoons. However, future research should relax this assumption and investigate the implications of the proposed solution on formation and stability of platoons.\\\\
\textbf{Assumption 3:}
The length of the schedule zone is sufficiently large so that a platoon can accelerate up to the speed limit and decelerate to complete stop. \\\\
We impose this assumption to ensure that the platoon entering the schedule zone with  speed less than the speed limit can reach the speed limit before it enters the merging zone of the intersection. It also ensures that the platoon entering the schedule zone with a speed equal to the speed limit will have the time to decelerate to a complete stop.
We assign the maximum speed for a platoon to enter the merging zone to be equal to the speed limit.
\subsection{Modeling Framework and Constraints}
Let $\displaystyle N(t) \in  \mathbb{N}$ be the number of platoons entering into the schedule zone at time $\displaystyle t \in \mathbb{R}^+$. The coordinator assigns a unique identification number $\displaystyle j \in \mathbb{N}$ to each platoon at the time they enter the schedule zone.  Let $\mathcal N(t)=\{1,\ldots N(t)\}$ be the queue of platoons inside the schedule zone. Let ${A}_j=\{1,\ldots,n_j\},~n_j\in \mathbb{N}$, be the number of vehicles in each platoon $\displaystyle j \in \mathcal N(t)$. We model each vehicle $\displaystyle i \in A_j$ as a double integrator,
\begin{align}
%\label{eqn}
\dot{p_i}=v_i(t), \nonumber \\ 
\dot{v_i}=u_i(t), 
\label{eq1}
\end{align}
where $\displaystyle p_i(t) \in \mathcal{P}_i$, $\displaystyle v_i(t) \in \mathcal{V}_i$,  $\displaystyle u_i(t) \in \mathcal{U}_i$ denote position, velocity, acceleration/deceleration. Let $\displaystyle x_i(t)=[p_i(t) \quad v_i(t)]^T$ denote the state of each vehicle $\displaystyle i \in A_j$. Let $\displaystyle t_{i}^{0}$ be the time at which vehicle $\displaystyle i \in A_j$ enters the schedule zone. Let ${x_i}^0=[p_i^0 \quad v_i^0]^T$ be the initial state where $\displaystyle  p_i^0=p_i(t_i^0)=0$, taking values in the state space $\displaystyle \mathcal{X}_i=\mathcal{P}_i \times \mathcal{V}_i$.  The control input and speed of each vehicle $\displaystyle i \in A_j$ is bounded with following constraints
\begin{gather}
%\label{eqn}
u_{\min} \leq u_i(t)  \leq u_{\max},  \label{eq2} \\
0 \leq v_{\min} \leq v_i(t) \leq v_{\max}, \label{eq3} 
\end{gather}
where $\displaystyle u_{\min}, u_{\max}$ are the minimum and maximum control inputs and $\displaystyle v_{\min}, v_{\max}$  are the minimum and maximum speed limits, respectively. 

\subsection{Modeling Left turn and right turns at an intersection}
\par
\begin{figure}[htbp!]
	\centering
	\includegraphics[width=3.2in]{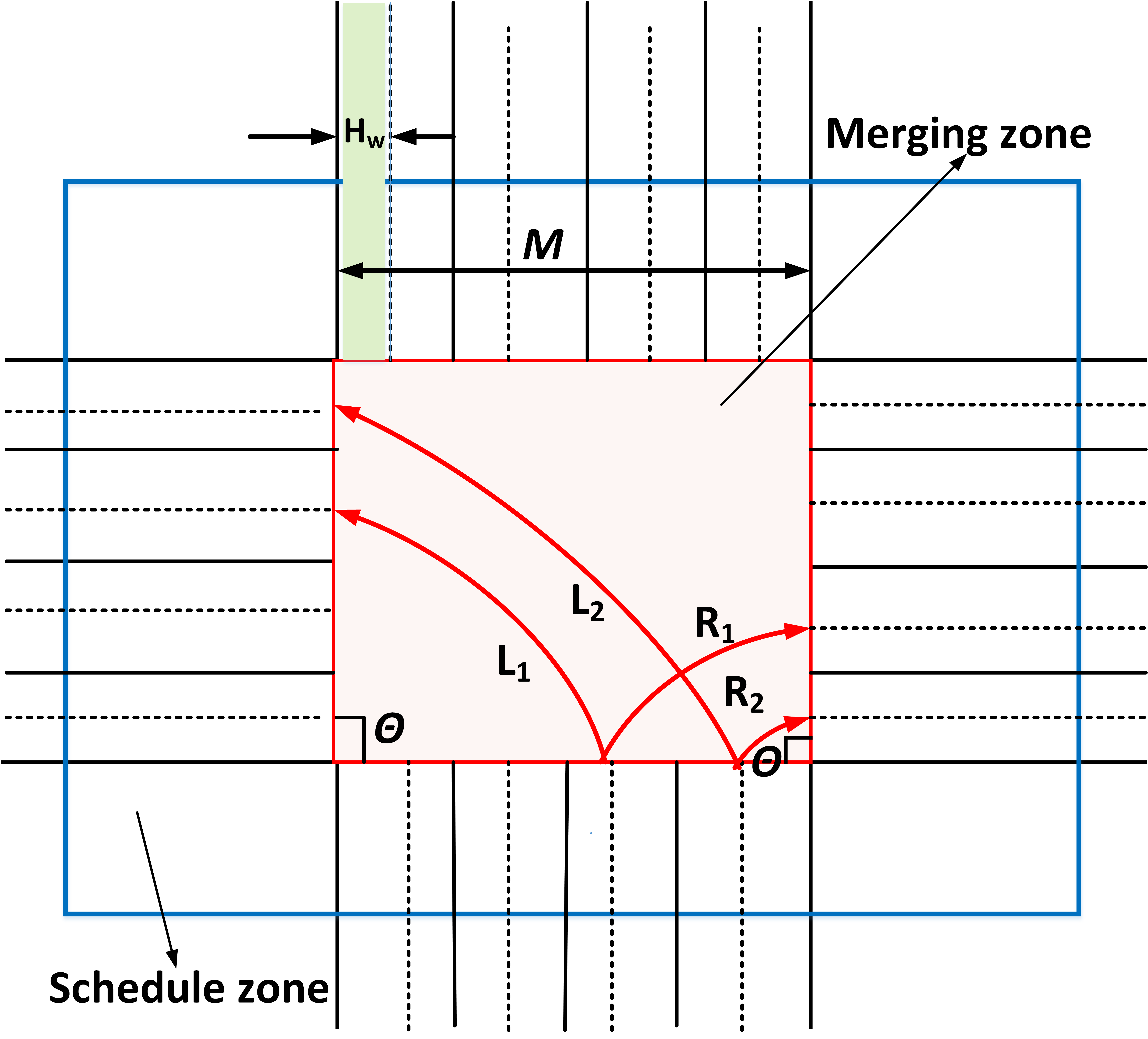}
	\caption{Intersection layout.}
	\label{arcL}
\end{figure}

 Let $D_j$ denote the routing decision (straight/left/right) of platoon $j$. Here, $D_j$ = S denotes the decision to go straight, $D_j$ = L denotes the decision to turn left,  and $D_j$ = R denotes the decision  to turn right at the intersection.
We consider an intersection layout as shown in Fig. \ref{arcL}.

The distance covered by a turning vehicle at the intersection is 
\begin{align}
d =\frac{\theta}{360}(2\pi{r}),
\end{align}
where $\theta$ is the angle subtended at the centre in radians and $r$ is the turning radius.

Let $W$ be number of the lanes in the road approaching the intersection. Let $H_w$ be half of a lane i.e.,
\begin{equation}
H_w=\frac{width~ of~ a~lane}{2}
\end{equation}

Let $H_r$ be the number of half lanes ($H_w$) from right end of the road and  $H_l$ be the number of half lanes ($H_w$) from left end of the road. The distance covered by a platoon $j$ that goes either straight, left or right inside merging zone is
\begin{equation}
d_j^*=\begin{cases}
 M, & \text{$D_j$ = S}, \\
1-\frac{H_r}{2W} \pi M, & \text{$D_j$ = L},\\
{1}-\frac{H_l}{2W} \pi M, & \text{$D_j$ = R}.
\end{cases}\label{eqdis} 
\end{equation}

In case of left turn $L_1$ and right turn $R_1$ at the  intersection (Fig. \ref{arcL}), $W=4$,  $H_l=5$ and $H_r=3$. The distance covered by $L_1$ and $R_1$ are
$\frac{5}{8}\pi M$ and $\frac{3}{8}\pi M$, respectively. In case of left turn $L_2$ and right turn $R_2$ at the  intersection (Fig. \ref{arcL}), $W=4$,  $H_l=7$ and $H_r=1$. The distance covered by $L_2$ and $R_2$ are $\frac{7}{8}\pi M$ and $\frac{1}{8}\pi M$, respectively.

The maximum allowable speed limit $v$ for   turning vehicles \cite{aashto2001policy} is 
\begin{align}
v=\sqrt{15R(0.1E+F)},
\end{align}
where $R$ is the effective centerline turning radius, E is the super-elevation (zero in urban conditions), and F is the side friction factor.
The maximum speed limit $v_{\max}^l$ and $v_{\max}^r$ of the platoons turning left  and right, respectively inside merging zone are 
\begin{align}
v_{\max}^l=\sqrt{15R_l(0.1E+F)},\\
v_{\max}^r=\sqrt{15R_r(0.1E+F)},
\end{align}
where $R_l$ and $R_r$ effective centerline turning radius of left turn and right turn at the intersection, respectively.
The maximum speed limit $v_{\max}$ of the platoon that goes either straight, turning left or right inside merging zone is denoted as, 

\begin{equation}
v_{\max}=\begin{cases}
 v_{\max}^s, & \text{$D_j$ = S}, \\
v_{\max}^l, & \text{$D_j$ = L},\\
v_{\max}^r, & \text{$D_j$ = R}.
\end{cases}\label{eqr} 
\end{equation}

\subsection{Upper-Level Optimal Framework for Coordination of platoons}

In this section, we discuss the upper-level optimization framework that yields the optimal schedule for the platoons to cross the merging zone with a minimum delay. The proposed framework is based on scheduling theory which addresses the allocation of jobs to the machines for a specified period of time aiming to optimize the performance measures.  A scheduling problem is described by the following notation $\mathcal{M}|\mathcal{C}|\mathcal{O}$, where $\mathcal{M}$ denotes machine environment, $\mathcal{C}$ denotes the constraints, and $\mathcal{O}$ denotes the objective function. We consider a job-shop scheduling problem where several jobs are processed in a single machine environment. Let $K\in \mathbb N$ be the number of jobs to be processed in a single machine. Let $t_{k}^p$  and $\displaystyle t_{k}^d$ be the processing time and deadline for each job $k \in K$. In a machine $\mathcal{M}$, if a job  $\displaystyle k$ starts at time  $\displaystyle t_{k}^s$ and completes at time $\displaystyle t_{k}^f$, then the completion time of job $k$ is $t_{k}^f=t_{k}^s+t_{k}^p$.\\\\
\textbf{Definition 1:}
 The lateness $\displaystyle \mathcal{L}_k$ of a job $k$ is defined as 
 \begin{align}
 \mathcal{L}_k \stackrel{\Delta}{=} t_{k}^f-t_{k}^d. \label{eq4}
 \end{align}
 The job-shop scheduling problem of minimizing maximum lateness in a single machine environment is represented as $\displaystyle 1||L_{\max}$ problem, where $1$ denotes a single machine and $L_{\max}$ denotes the maximum lateness. A schedule is said to optimal if it minimizes $\displaystyle \max_{k} \mathcal{L}_k$, i.e., the maximum lateness of jobs.
\par In our proposed framework, we model the intersection as a single machine and the platoons as jobs. Based on Assumption 1, the vehicles form stable platoon and each stable platoon is considered as a job. The processing time is the time taken by the job to be completed in a machine. We model the processing time of a job as \emph{passing time} of platoons, i.e., the time taken by the platoons at the maximum speed to exit the intersection. The deadline of a job is the time before which it must be completed in a machine.  We model the deadline of the job as \emph{deadline} of the platoons, i.e., the time taken by the platoons at their initial speed to exit the intersection. Then, each platoon leader solves $1||L_{\max}$ scheduling problem in a single machine environment (intersection) to find the optimal schedule to enter the merging zone of the intersection that minimizes maximum lateness, i.e., travel delay. \\\\
\textbf{Definition 2:}
 Let $\displaystyle t_{j}^{0}$ and  $\displaystyle t_{j}^{m}$ be the time at which the platoon $\displaystyle j \in \mathcal{N}(t)$ enters the schedule zone and merging zone, respectively.  The arrival time period $t_{j}^{a}$  of the platoon $\displaystyle j\in \mathcal{N}(t)$ at the merging zone  is  
\begin{align}
t_{j}^a \stackrel{\Delta}{=} t_{j}^{m}-t_{j}^{0}.
\label{eq5} 
\end{align}
\textbf{Definition 3:}
Let  $\displaystyle t_{j}^{e}$ be the time at which the platoon  $\displaystyle j \in \mathcal{N}(t)$ exits the merging zone. The crossing time period $\displaystyle t_{j}^{c}$ of a platoon $\displaystyle j \in \mathcal{N}(t)$ is 
\begin{align}
t_{j}^{c}\stackrel{\Delta}{=}t_{j}^{e}-t_{j}^{m}.
\label{eq6} 
\end{align} 
\textbf{Definition 4:}
 The passing time $\displaystyle t_{j}^p$ of a platoon $\displaystyle j \in \mathcal{N}(t)$ at the intersection is 
\begin{align}
t_{j}^p \stackrel{\Delta}{=}  t_{j}^{a} + t_{j}^{c}.
\label{eq7} 
\end{align}
\par
\begin{figure*}[htbp!] 
 	\centering
 \includegraphics[width=5in, height=4in]{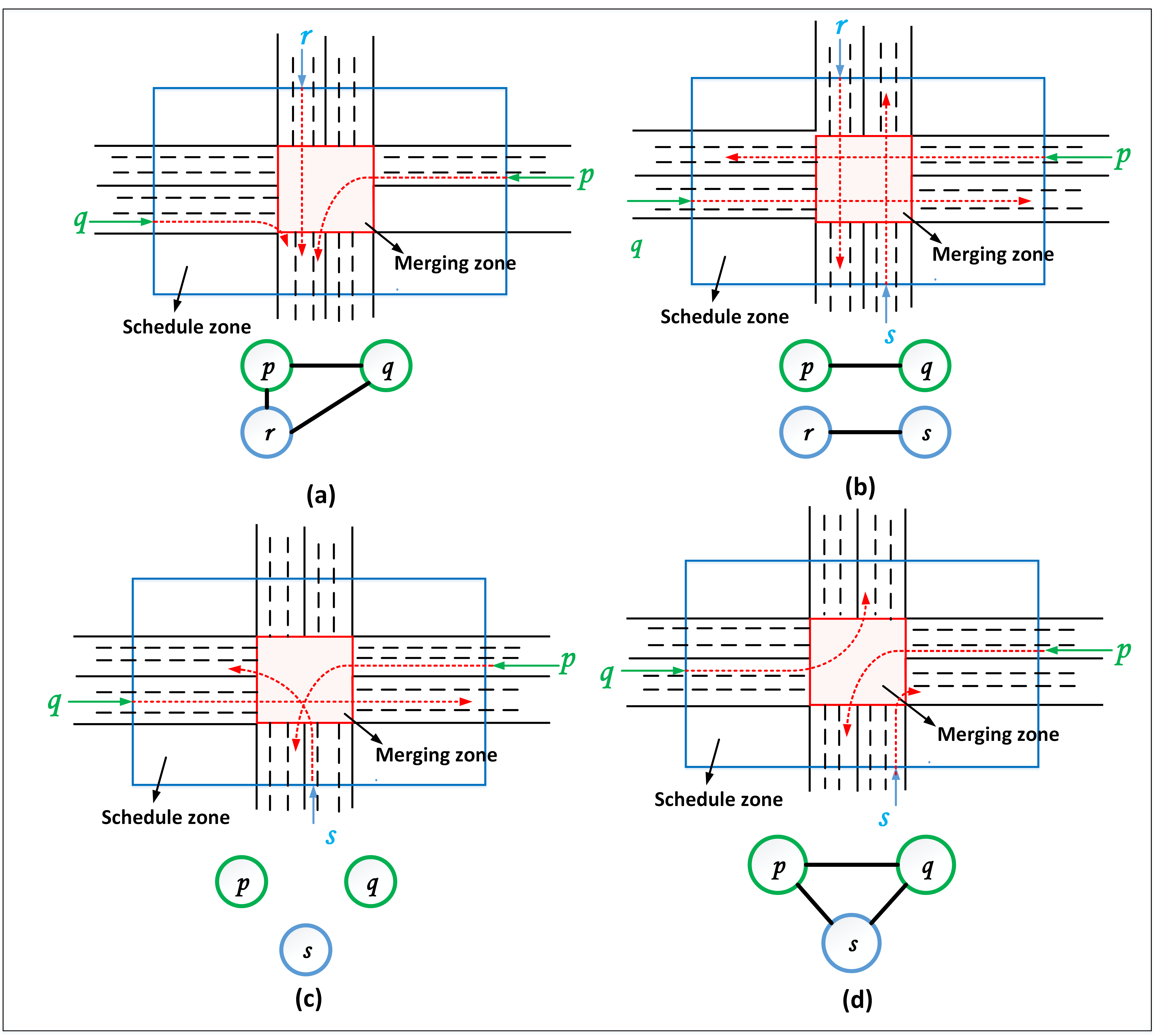}
 	\caption{Intersection with traffic movements.}
 	\label{intersection}
 \end{figure*}

 We consider two cases for computing the passing time of platoons at the time they enter the schedule zone. In Case 1, the platoon enters the schedule zone while cruising with the speed limit. In Case 2,  the platoon enters the schedule zone with speed that is less than the speed limit.  Let $\displaystyle v_{j}^0= v_j(t_j^0)$ be the initial speed of the platoon $\displaystyle j \in \mathcal{N}(t)$, i.e., speed at which the platoon enters the schedule zone.  Let $|A_j|$ be the cardinality of $A_j$.  Let $\displaystyle t_{j}^h$ be the time headway between the vehicles $\displaystyle i$ and $\displaystyle (i-1) \in A_j$ in the platoon. For instance, if there is a platoon of $4$ vehicles and all moving at constant speed of $10~m/s$ with uniform time headway of $1.2~sec$, the space headway between vehicles in a platoon is $(10\times1.2)=12~m$. Let $\displaystyle t_{c}$ be the clearance time interval, i.e., a safe time gap provided between exit and entry of platoons at the merging zone to ensure safety of platoons.\\\\
\textbf{Case 1: $v_{j}^0= v_{\max}$}\\
Using  $\displaystyle v_{\max}$, we compute the arrival time    
\begin{gather}
t_{j}^{a*}=\frac{S}{v_{\max}},\label{eq8}
\end{gather}
and the crossing time of platoons 
\begin{gather}
t_{j}^{c*}=\frac{d^*}{v_{\max}} + (|A_j|-1) \times t_{j}^h+t_{c}.\label{eq9}
\end{gather}\\\\
\textbf{Case 2: $ v_{j}^0 < v_{\max}$} \\ 
We compute $\displaystyle {t_{j}^{a}}^*$ using the time taken by the platoon to accelerate to the speed limit applying its maximum acceleration.
Let $\displaystyle t_{j}^s$ be the time taken by the platoon $j$ to accelerate to the speed limit, then we have 
\begin{align}
t_{j}^s=\frac{v_{\max}-v_{j}^0}{u_{\max}}\label{eq10}.
\end{align}
 Let  $\displaystyle d_{j}^s$ be the distance traveled during acceleration,  then we have
 \begin{align}
d_{j}^s=\frac{{(v_{\max})^2}-(v_{j}^0)^2}{2u_{\max}}. \label{eq11}
 \end{align}
 Based on Assumption 3, the platoons will reach the speed limit at time $t \leq {t_{j}^{a}}^*$  and 
 \begin{gather}
{t_{j}^{a}}^*=t_{j}^s+\frac{S- d_{j}^s}{v_{\max}},\label{eq12}
\end{gather}
and $t_{j}^{c*}$ is computed using (\ref{eq9}).\\\\
\textbf{Definition 5:}
Let $\displaystyle t_{j}^{a'} $ be the time taken by the platoon $\displaystyle j$ to reach the the merging zone while cruising with their initial speed. The deadline $\displaystyle t_{j}^d$  of the platoon $\displaystyle j$ to completely cross the intersection is defined as
\begin{align}
t_{j}^d \stackrel{\Delta}{=}  t_{j}^{a'} + t_{j}^{c*}.
\label{eq14}
\end{align}
We compute $t_{j}^{a'}$  as 
\begin{gather}
t_{j}^{a'}=\frac{S}{v_{j}^0},\label{eq15}
\end{gather}
and the crossing time $\displaystyle t_{j}^{c*}$ of the platoon using (\ref{eq9}).\\\\
\textbf{Definition 6:}
Let $\Gamma_i$ and $\Gamma_j$ be the path of platoons $i$ and $j \in \mathcal{N}(t)$, respectively. The platoons $i$ and $j$ are said to be compatible if $\Gamma_i \cap \Gamma_j= \varnothing$,  i.e., paths of platoons $i$ and $j$ are non-conflicting and can be given right-of-way concurrently inside the merging zone.\\
The compatibility between the paths of the platoons can be modeled as a compatibility graph.\\\\
\textbf{Definition 7:}
A compatibility graph $\mathcal{G}_c=(\mathcal{V, E})$ is an undirected graph where $\displaystyle \mathcal V$ is the set of vertices and  $\mathcal{E}$ is the set of edges. The adjacency matrix $A=[a_{ij}]$ of compatibility graph  $\mathcal{G}_c$ can be defined as
\begin{equation}
a_{ij}=\begin{cases}
 1, & \text{if paths of platoons}~i ~\text{and} ~j~ \text{do not conflict}, \\
0, & \text{if paths of platoons}~i ~\text{and} ~j~ \text{conflict}.
\end{cases}\label{eqp} 
\end{equation}
For example, we consider an intersection with traffic movements as shown in Fig. \ref{intersection}. 
 In Fig. 3a, let  $p, q$, and $r \in \mathcal{N}(t)$ be the platoons entering the schedule zone at time $t$. 
Here, $\mathcal V=\{p,q,r\}$ is the vertex set of compatibility graph $\displaystyle \mathcal{G}_c$. An edge $e$ connects the vertices \{$p, q$\}, \{$q, r$\}, and \{$p, r$\} since their paths are non-conflicting inside the merging zone. \\\\
 \textbf{Definition 8:} A clique $\mathcal{C}$ of  $\mathcal{G}_c$ is a subset of the vertices $\mathcal{C} \subseteq \mathcal{V}$ such that each vertex in $\mathcal{C}$ is adjacent to all other vertices in $\mathcal{C}$. \\\\
\textbf{Definition 9:} A maximal clique $\mathcal{M}$ is the clique that consists of a set of vertices $\mathcal{M} \subseteq \mathcal{V}$ which is not a subset of any other cliques in the undirected graph $\mathcal{G}_c$.\\\\
In Fig. 3a, the set of vertices \{$p$, $q$, $r$\} forms the maximal clique of the compatibility graph $\mathcal{G}_c$. In our framework, the maximal cliques represents the groups of compatible platoons that can be given right-of-way concurrently inside the merging zone. Therefore, there is one group of compatible platoons. In Fig. 3b, The set of vertices \{$p$, $q$\} and \{$r$, $s$\} are the maximal cliques of the compatibility graph $\mathcal{G}_c$. Therefore, there are two groups of compatible platoons. In Fig. 3c, The set of vertices \{$p$\}, \{$q$\}, and \{$s$\} are the maximal cliques of the compatibility graph $\mathcal{G}_c$. Therefore, there are three groups of compatible platoon. In Fig. 3d, The set of vertices \{$p$, $q$, $s$\} form the  the maximal clique of the compatibility graph $\mathcal{G}_c$. Therefore, there is one group of compatible platoons. We generate group of compatible platoons as illustrated in Fig. 3 for various combinations of routing decisions of platoons. \\\\
\textbf{Definition 10:} Let $G=\{1,\ldots, n\}$ be the set of groups of compatible platoons. The completion time $\displaystyle t_{G}^f$ is defined as the time taken by all groups of platoons in $\displaystyle G$ to completely exit the merging zone of the intersection.\\\\
 In the upper-level optimization framework, we model the problem of coordinating platoons at the intersection as a job-shop scheduling problem, the solution of which yields the optimal schedule for each platoon to cross the intersection through four algorithms. The proposed framework uses the information about passing time and deadline of each platoon entering the schedule zone. Algorithm $1$ computes the passing time of each platoon given the speed limit, geometric information of the intersection, and attributes ($n_j$, $N_{link}$, $N_{lane}$, $D_j$, $p_j$, $v_j$) of the platoon $j \in \mathcal N(t)$. Algorithm $2$ computes the deadline of each platoon to cross the intersection. Next, Algorithm $3$ categorizes the platoons into groups of compatible platoons and computes the passing time, deadline, and crossing time of each group. The earliest due date principle \cite{jackson1955scheduling} for scheduling jobs in a single machine environment is optimal in minimizing the maximum lateness of the jobs. We adapt the earliest due date principle in Algorithm 3 to find an optimal sequence of platoons to reduce the delay. Finally, Algorithm 4 computes the time of entry for each platoon inside the merging zone. The upper-level optimal framework thus yields an optimal schedule to reduce the delay of platoons which is equivalent to  $1||L_{\max}$ scheduling problem and a proof is presented below.\\\\
\textbf{Theorem 1:} \emph{The schedule $\displaystyle S$ in the non-decreasing order of  deadline $\displaystyle t_{g}^d$ of each group $g \in G$ is optimal in minimizing the travel delay, i.e., the maximum lateness of the platoons}.\\\\
\textbf{Proof:}
Let $t_{g}^d$ be deadline of the each group $g\in G$ and the groups are arranged in non-decreasing order of their deadlines.  Let consider  two groups of platoons $\displaystyle i$ and $\displaystyle j$ arriving at the schedule zone at time $\displaystyle t$. Let the lateness of the platoon group $\displaystyle i$ and $\displaystyle j$ be $\displaystyle \mathcal{L}_{i}$ and  $\displaystyle \mathcal{L}_{j}$, respectively. Suppose there exists a schedule $\displaystyle S'$ in which the group $\displaystyle i$  enters  the merging zone before the  group $\displaystyle j$ and  $\displaystyle t_{i}^d > t_{j}^{d}$. Then,
\begin{gather}
 \mathcal{L}_{i}^{S'}=(t+t_{i}^{p})-t_{i}^{d},\label{eq17}\\
\mathcal{L}_{j}^{S'}=(t+t_{i}^{p}+t_{j}^{p})-t_{j}^{d},\label{eq18}
\end{gather}
which implies,
\begin{gather}
\quad  \mathcal{L}_{j}^{S'}> \mathcal{L}_{i}^{S'}. \quad \label{eq19}   
\end{gather}
Suppose there is an another schedule $\displaystyle S$, in which group $\displaystyle j$ enters the merging zone before the group $\displaystyle i$ at time $\displaystyle t$. Then we have
\begin{gather}
\mathcal{L}_{j}^S=(t+t_{j}^{p})-t_{j}^{d} < \mathcal{L}_{j}^{S'}, \label{eq20}\\
\mathcal{L}_{i}^S=(t+t_{j}^{p}+t_{i}^{p})-t_{i}^{d} < \mathcal{L}_{j}^{S'}.\label{eq21} 
\end{gather}
Thus,
\begin{gather}
\quad \max\{\mathcal{L}_{j}^S,\quad\mathcal{L}_{i}^S\}  \leq  \max\{\mathcal{L}_{i}^{S'},\quad\mathcal{L}_{j}^{S'}\}. \label{eq22}  
\end{gather}
Thus, the schedule $S$ is optimal if and only if the deadlines $t_{g}^d$ of groups of compatible platoons have been sorted in a non-decreasing order.
\hspace*{86mm}\qedsymbol
\begin{algorithm}
	\caption{Compute passing time of platoons}
	 \hspace*{1mm} \textbf{Input:}  $v_{j}^0$, $|A_j|$, $t_{j}^h$ of each platoon $j$, $v_{\max}$, $t_{c}$, $S$,  $M$, \hspace*{12mm} $u_{\max}$.\\
\hbox{\textbf{Output:} $t_{j}^{a*}$,  $t_{j}^{c*}$ and $t_{j}^p$  of each platoon $j$.}\\
%\vspace{2mm}
  \{$\triangleright$ Computation of arrival time\}
	\begin{algorithmic}[1]
	    % \Comment{Computation o  f arrival time}
	    % 
	   % \vspace{2mm}
	   % \Statex \ $\triangleright$ Computation of arrival time\ 
	    \FOR {$j=1~ to~ N$}
	%	 \vspace{2mm}
		\IF {$v_{\max}=v_{j}^0$}
	%	\vspace{2mm}
		\STATE $t_{j}^{a*} \gets S/v_{\max}$ %\COMMENT{Computation of arrival time} 
	%	\vspace{2mm}
		\ELSIF {$v_{\max}<v_{j}^0$}
	%	\vspace{2mm}
		\STATE $t_{j}^{s} \gets (v_{\max}-v_{j}^0)/ u_{\max}$
	%	\vspace{2mm}	
		\STATE $d_{j}^s \gets [(v_{\max}^2)-(v_{j}^0)^2]/ 2u_{\max}$
	%	\vspace{2mm}
		\STATE $t_{j}^{a*} \gets t_{j}^s+(S-d_{j}^s)/v_{\max}$
	%	\vspace{2mm} 
		\ENDIF		    
		\ENDFOR\\
	%	\vspace{2mm}
		\COMMENT{$\triangleright$ Computation of crossing time}
	%	 \vspace{2mm}
		\FOR {$j=1~to~N$}
	%	\vspace{2mm}
		\STATE $t_{j}^{c*} \gets d^*/v_{\max}+(|A_j|-1)*t_j^{h}+t_{c}$
	%	\vspace{2mm}
		\ENDFOR\\
	%	\vspace{2mm}
\COMMENT{$\triangleright$ Computation of passing time}
		%\vspace{2mm}
		\FOR {$j=1~to~N$}
	%	\vspace{2mm}
		\STATE  $t_{j}^p \gets t_{j}^{a*} +t_{j}^{c*}$
	%	\vspace{2mm}
		\ENDFOR
	\end{algorithmic} 
\end{algorithm}	

\begin{algorithm}
	\caption{Compute deadline for platoons to exit the merging zone}	
	 \hspace*{1mm} \textbf{Input:}  $v_{j}^0$, $t_{j}^{c*}$ of each platoon $j$, $S$.\\
	\hbox{\textbf{Output:} {$t_{j}^d$ of each platoon $j$.}}\\
%	\hspace*{6mm} 
\{$\triangleright$ Computation of arrival time\}
	\begin{algorithmic}[1]
	%\vspace{2mm}
	  %\Statex \ $\triangleright$ Computation of arrival time\
	  %	\vspace{2mm}
		\FOR {$j=1~ to~ N$}
	%	\vspace{2mm}
		\STATE  $t_{j}^{a'} \gets S/v_{j}^0$
		%\vspace{2mm}
		\ENDFOR\\
	%	\vspace{2mm}
	 % \Statex \$\triangleright$ Computation of deadline\
	 \COMMENT{$\triangleright$ Computation of deadline}
	  %	\vspace{2mm}
		\FOR {$j=1~ to~ N$}
	%	\vspace{2mm}
		\STATE  $t_{j}^d \gets t_{j}^{a'}+ t_{j}^{c*}$
		%\vspace{2mm} 	 
		\ENDFOR 
	\end{algorithmic} 
\end{algorithm}

\begin{algorithm}
	\caption{Compute the groups of compatible platoons and optimal sequence}	
	\hspace*{1mm} \textbf{Input:} compatibility graph $\mathcal {G}_c$, $t_{j}^p$, $t_{j}^{c*}$, $t_{j}^d$ of each Platoon $j$.\\
	\hbox{\textbf{Output:} {groups of compatible platoons $G$, $t_{g}^d$, $t_{g}^p$, $t_{g}^{c*}$ of }}\\ \hbox{ \hspace*{11mm} each group $g$ and optimal sequence.}
	\vspace{-4mm}
	\begin{algorithmic}[1]
      \STATE  $G \gets$ maximal cliques of compatibility graph, $\mathcal {G}_c$.\\
      	%\vspace{2mm}
	  \COMMENT {$\triangleright$ Computation of passing time, crossing time and deadline of each group of platoon}
	  %	\vspace{2mm}
     \STATE Initialize Array $deadline \gets 0$
     %	\vspace{2mm}
     \STATE Initialize Array $passingTime \gets 0$
     %	\vspace{2mm}
     \STATE Initialize Array $exitTime \gets 0$
     %	\vspace{2mm}
     \STATE Initialize variable $i \gets 1$
     %	\vspace{2mm}
      \FOR {each $g~ in ~G$}    
      %	\vspace{2mm}
       	\FOR {$j=1~ to~ N$}
       	
       	%	\vspace{2mm}
      	\IF {$j \in g$}
      	%	\vspace{2mm}
        \STATE	$deadline[i] \gets t_{j}^d$
       % \vspace{2mm}
        \STATE $passingTime[i] \gets t_{j}^{p}$
       % \vspace{2mm}
        \STATE $exitTime[i] \gets t_{j}^{c*}$
        %	\vspace{2mm}
        \STATE $i\gets i+1$
       % 	\vspace{2mm}
      	\ENDIF
      	%	\vspace{2mm}
      	\ENDFOR
      	%	\vspace{2mm}
      	\STATE $t_{g}^d \gets max(deadline)$
      %	\vspace{2mm}
      	\STATE $t_{g}^{c*} \gets max(exitTime)$
      %	\vspace{2mm}
      	\STATE $t_{g}^{p} \gets max(passingTime)$
      	%	\vspace{2mm}
      	\STATE $i\gets 1$
      	%	\vspace{2mm}
      	\ENDFOR	
      	%	\vspace{2mm}
      \STATE optimalSequence = Platoon groups sorted in non-decreasing order of deadline $\displaystyle t_{g}^d$ 
	\end{algorithmic} 
\end{algorithm}
\
\begin{algorithm}
	\caption{Compute time of entry for each platoon}	
	\hspace*{1mm} \textbf{Input:} current time $t$, optimalSequence, $t_{g}^{p}$, $t_{g}^{c*}$, $t_{g}^d$ of each group $g$, $t_{j}^{a*}$ of each platoon $j$.\\
	\hbox{\textbf{Output:} {$t_{j}^{m}$ of each platoon $j$.}} \\
 \hspace*{4mm}	\{$\triangleright$ Computation of time of entry for first group in the optimal sequence\}
	\begin{algorithmic}[1]
	%\vspace{2mm}
	    % \Statex \ $\triangleright$ Computation of time of entry for first group in the optimal sequence\
	  %	\vspace{2mm}
		\STATE Initialize variable $t_{G}^f \gets 0$
	%	\vspace{2mm}
		\STATE $l=optimalSequence[1]$
	%	\vspace{2mm}
		\FOR {each $j \in l$}
	%	\vspace{2mm}
		\STATE $t_{j}^{m} \gets t + t_{j}^{a*}$
	%	\vspace{2mm}
		\ENDFOR
		%\vspace{2mm}
		\STATE  $t_{l}^e\gets t+ t_{l}^{p}$		%\vspace{2mm}
		\IF {$t_{j}^{m} \leq t_{G}^f$}
	%	\vspace{2mm}
		\FOR {each $j \in l$}		
	%	\vspace{2mm}
		\STATE	$t_{j}^{m} \gets t_{G}^f$
	%	\vspace{2mm}
		\ENDFOR
	%	\vspace{2mm}
		\STATE $t_{l}^e \gets t_{G}^f + t_{l}^{c*}$
	%	\vspace{2mm}
		\ENDIF\\
	%	\vspace{2mm}	
		\COMMENT {$\triangleright$ Computation of time of entry for groups $2$ to $N$ in the optimal sequence}
	%	\vspace{2mm}
		\FOR {$k=2~ to~ N$}
	%	\vspace{2mm}
		\STATE $r=optimalSequence[k]$
	%	\vspace{2mm}
		\FOR {each $j \in r$}
	%	\vspace{2mm}
		\STATE	$t_{j}^{m} \gets t_{r-1}^{e}$
	%	\vspace{2mm}
		\ENDFOR
	%	\vspace{2mm}
		\STATE $t_{r}^e\gets t+ t_{r}^{c*}$
	%	\vspace{2mm}
		\STATE $lastG=optimalsequence[N]$
	%	\vspace{2mm}
		\ENDFOR
	%	\vspace{2mm}
		\STATE $t_{G}^f \gets t_{lastG}^{e}$ 
	\end{algorithmic} 
\end{algorithm}

The platoon leader runs the algorithm at every time step after entering the schedule zone. The platoon leaders inside the schedule zone communicate with each other and find the time they can enter the merging zone. Then, they derive their optimal control input to enter the merging zone.  When a new platoon enters the schedule zone,  the platoon leader along with the other platoon leaders run the algorithm to find the new schedule excluding the platoons that have entered the merging zone. After entering the merging zone, the platoon leader does not communicate to the other platoon leaders. The new platoon in the schedule zone  computes its time of entry into the merging zone based on the platoons inside the schedule zone and the crossing time of the platoons that have entered the merging zone.  The platoons inside the schedule zone should not enter the merging zone until the crossing platoon has completely exited the merging zone.  The new platoon leaders inside the schedule zone fetch the crossing time of platoons that have entered the merging zone from the coordinator which is acting as a database. 

\subsection{Low-Level Framework for Optimal Control Problem}
In the low-level optimization framework, each platoon derives the optimal control input based on the time of entry into the merging zone designated by the upper-level framework. The platoons enter the merging zone at the designated time  and in the order of the optimal sequence provided by upper-level framework. The platoons are allowed to cross the intersection based on their position in the optimal sequence. In that case, some platoons enter the merging zone at its earliest arrival time. Other platoons inside the schedule zone  have to wait for the other platoons inside the merging zone to exit the merging zone. The platoons waiting for other platoons to exit the merging zone derive energy optimal trajectory to enter the merging zone at the time specified by the upper-level optimal framework. If the time that a platoon enters inside the merging zone is equal to its earliest arrival time, then the leader derives the optimal control input by solving the time optimal control problem. If the time of entry of a platoon inside the merging zone is greater than the earliest arrival time, then the  leader derives the optimal control input by solving an energy optimal control problem.  After deriving the optimal control input, based on Assumption 2, each platoon leader communicates the time of entry and optimal control input (acceleration/deceleration) to the followers until the last vehicle in the platoon exit the merging zone.
  \subsubsection{Time Optimal Control Problem}
For each leader  $\displaystyle l \in A_j$ of the platoon $j \in \mathcal{N}(t)$, we define the following optimal control problem 
\begin{align}
\min_{u_{l}\in \mathcal{U}_l}J_1(u_l(t))=\int_{t_l^0}^{t_l^m}dt= t_{l}^{m}-t_{l}^{0},\label{eq23}
\end{align}
subject to: (\ref{eq1}), (\ref{eq2}), (\ref{eq3}), $\displaystyle p_l(t_l^0)=0$, $\displaystyle p_l(t_{l}^{m})=L_{s}$, \\ and given  $\displaystyle t_l^0, v_{l}^0, t_{l}^{m}$,\\\\ where $t_l^0$ is time that the platoon leader enters the schedule zone, and $t_{l}^{m}$ is time that the platoon leader enters the merging zone.

\subsubsection{Energy Optimal Control Problem}
For each leader  $\displaystyle l \in A_j$ of the platoon $j \in \mathcal{N}(t)$, we define the following optimal control problem 
\begin{align}
\min_{u_{l}\in \mathcal{U}_l} J_2(u_l(t))=\frac{1}{2}\int_{t_l^0}^{t_l^m} u_l^2(t) \; dt, \label{eq24} 
\end{align}
subject to: (\ref{eq1}), (\ref{eq2}), (\ref{eq3}),  $\displaystyle p_l(t_l^0)=0$,  $\displaystyle p_l(t_{l}^{m})=L_{s}$ and given  $\displaystyle t_l^0, v_{l}^0, t_{l}^{m}$,
\\\\
 where $t_l^0$ time that platoon leader enters the schedule zone, and $t_{l}^{m}$ is time that the platoon leader enters the merging zone. 

\section{Analytical Solution}
In this section, we derive the closed-form analytical solutions for the time and energy optimal control problems for each platoon leader $\displaystyle l \in A_j$.

\subsection{Analytical Solution of the Time Optimal control problem}
We apply Hamiltonian analysis for deriving the analytical solution of the time optimal control problem.  For each leader $l\in A_j$, the Hamiltonian function with the state and control constraints is
 \begin{gather}
H_l(t,p_l(t),v_l(t),u_l(t))=1 + {\lambda_l^p}{ v_l(t)}+ \lambda_l^v u_l(t) \nonumber \\
+ \mu_l^a (u_l(t)-u_{\max}) + \mu_l^b (u_{\min}-u_l(t))\nonumber \\+\mu_l^c (v_l(t)-v_{\max}) + \mu_l^d (v_{\min}-v_l(t)), \label{eq25} 
\end{gather}
where $\displaystyle {\lambda_l^p}$ and $\displaystyle {\lambda_l^v}$ are costates and  $\displaystyle {\mu_l^a}$,  $\displaystyle {\mu_l^b}$, $\displaystyle {\mu_l^c}$, and $\displaystyle {\mu_l^d}$ are lagrange multipliers.

We consider that the state constraint (\ref{eq3}) is not active, and therefore $\mu_l^c$ and $\mu_l^d=0$. Then,
\begin{align}
\lambda_l^p &= a_l, \label{eq32} \\
\lambda_l^v &= -a_lt+b. \label{eq33} 
\end{align}
From \eqref{eq25}, the optimal control input is,
\begin{equation}
u_l^*(t)=\begin{cases}
u_{\min}, & \lambda_l^v>0,\\
u_{\max}, & \lambda_l^v<0.
\end{cases}\label{eq34}
\end{equation}

 We consider two cases while platoons are entering the schedule zone of the intersection.\\
\textbf{Case 1:} If the platoon enters the schedule zone with $v_l(t)\le v_{\max}$, then from \eqref{eq34} we have
\begin{equation}
u_l^*(t)=\begin{cases}
u_{\max}, & \text{if $t_{l}^{0} \leq t \leq t_{l}^{0}+ t_{l}^s$},\\
0, & \text{if $t_{l}^{0}+ t_{l}^s \leq t \leq t_{l}^{m}$}.
\end{cases}\label{eq41} 
\end{equation}
Substituting (\ref{eq41}) in (\ref{eq1}), we can compute the
optimal position and velocity, 
\begin{align}
p_l^*(t) &=\frac{1}{2}u_lt^2+b_lt+c_l, \label{eq42}\\
v_l^*(t) &=u_lt+b_l,\quad \label{eq43} t \in [t_{l}^{0},~  t_{l}^{0}+ t_{l}^s],\\
p_l^*(t) &=v_{\max}t+d_l,\\
v_l^*(t) &=v_{\max},\quad \label{eq44}  t \in [t_{l}^{0}+ t_{l}^s, \quad t_{l}^{m}].
\end{align}
where $ b_l$, $c_l$, and $d_l$ are integration constants. We can compute these constants using the initial and final conditions in (\ref{eq23}).\\
\textbf{Case 2:} 
 If the platoon enters the schedule zone with $v_l(t)= v_{\max}$, then from \eqref{eq34} we have 
\begin{equation}
u_l^*(t)= 0,\quad{ t \in [ t_{l}^{0},~ t_{l}^{m}]}.
\label{eq46} 
\end{equation}
Substituting (\ref{eq46}) in (\ref{eq1}), we can compute the optimal position and velocity,
\begin{align}
p_l^*(t) &= v_{\max}t+d_l,\label{eq47}\\
v_l^*(t) &= v_{\max},~~~\label{eq48} t \in [t_{l}^{0},~ t_{l}^{m}].
\end{align}
where $ d_l$ is integration constant. We can compute the constant using the initial and final conditions in (\ref{eq23}). The complete solution of the optimal control problem with state and control constraints is presented in \cite{chalaki2019hierarchical}. 
  
\subsection{Analytical Solution of the Energy Optimal Control Problem}
For the analytical solution of energy optimal control problem, we apply Hamiltonian analysis with inactive state and control constraints.  We formulate the Hamiltonian function for each platoon leader $l\in A(t)$ as follows
\begin{gather}
H_l(t,p_l(t),v_l(t),u_l(t))=\frac{1}{2} u_l^2(t) + {\lambda_l^p}{ v_l(t)}+ \lambda_l^v u_l(t) \nonumber \\
+ \mu_l^a (u_l(t)-u_{l,\max}) + \mu_l^b (u_{l, \min}-u_l(t))\nonumber \\+\mu_l^c (v_l(t)-v_{l, \max}) 
+ \mu_l^d (v_{l, \min}-v_l(t)),\label{eq49}
\end{gather}
where $\displaystyle {\lambda_l^p}$ and $\displaystyle {\lambda_l^v}$ are costates, and  $\displaystyle {\mu_l^a}$,  $\displaystyle {\mu_l^b}$, $\displaystyle {\mu_l^c}$, and $\displaystyle {\mu_l^d}$ are the lagrange multipliers. Since the control and state constraints are not active, $\mu_l^a = \mu_l^b = \mu_l^c = \mu_l^d = 0$. The optimal control input based on  \cite{malikopoulos2018decentralized} will be
\begin{gather}
u_l^*=a_lt+b_l,~ \label{eq50}t \in [t_{l}^{0},~ t_{l}^{m}].
\end{gather}
Substituting (\ref{eq50}) in (\ref{eq1}), we can find the optimal position and velocity,
\begin{align}
p_l^* &= \frac{1}{6}a_lt^3+\frac{1}{2}b_lt^2+c_lt+d_l,~t \in [t_{l}^{0},~~ t_{l}^{m}], \label{eq51}\\
v_l^* &= \frac{1}{2}a_lt^2+b_lt+c_l, ~~t \in [t_{l}^{0},~ t_{l}^{m}],\label{eq52}
\end{align}
where $a_l,~ b_l,~ c_l$, and $d_l$ are integration constants. We can compute these constants using initial and final conditions, i.e., $p_l(t_{l}^{0})=p_l^0, v_l(t_{l}^{0})=v_l^0$, and $p_l(t_{l}^{m})=p_l^m,  v_l(t_{l}^{m})=v_l^m$. 
 
\section{Simulation Framework and Results}
\begin{figure}[htbp!]
 	\centering
 	\includegraphics[width=3in]{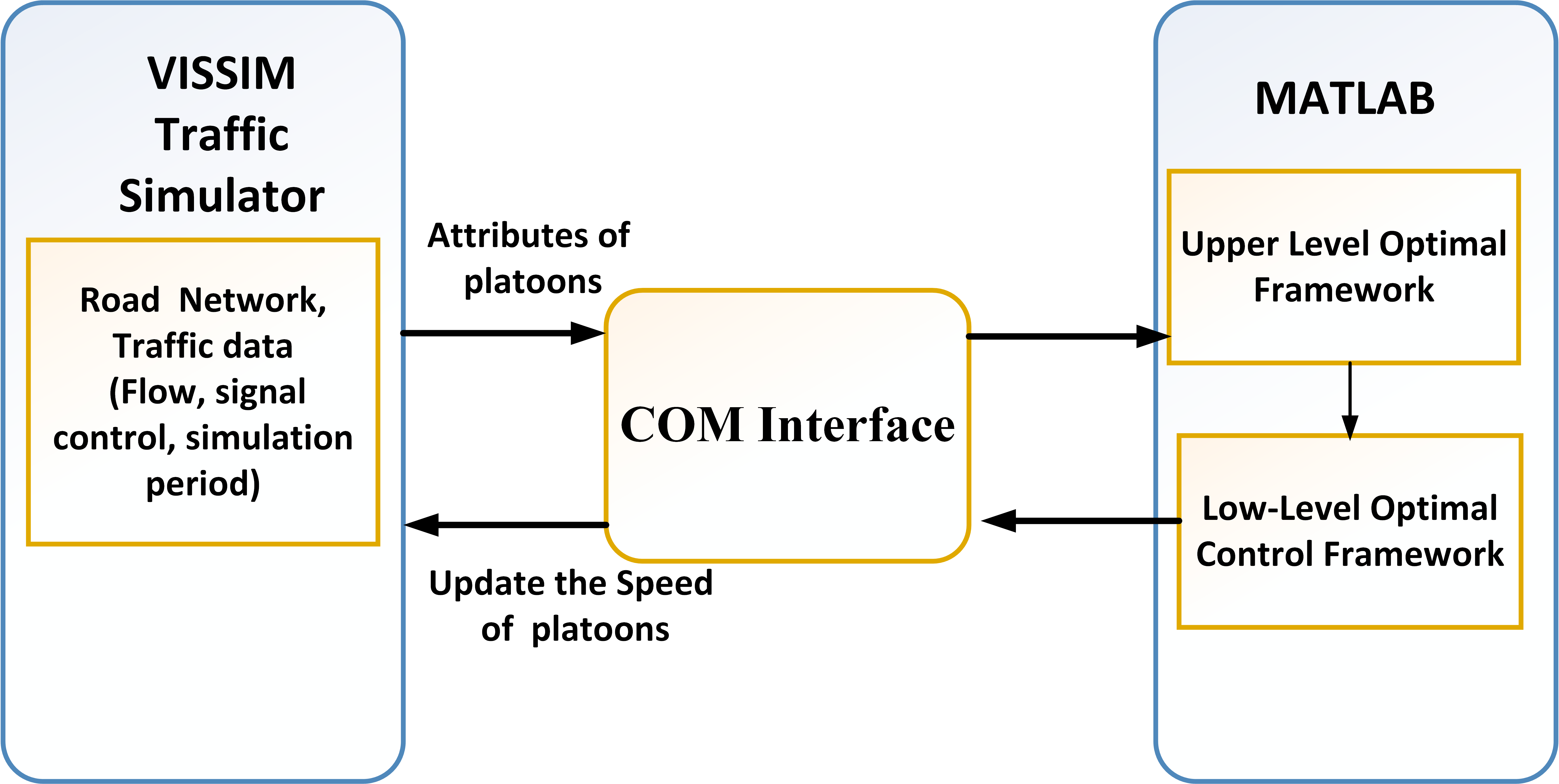}
 	\caption{Simulation framework in VISSIM-MATLAB environment.}
 	\label{simframe}	
 \end{figure}

We present the simulation framework (Fig. \ref{simframe}) using a VISSIM-MATLAB environment. In our simulation study, we model the intersection using VISSIM 11.00 traffic simulator \cite{PTV}. The length of the schedule zone $S$ and the merging zone $M$ for the intersection is $\displaystyle 200~m$ and $\displaystyle 50~m$, respectively. We designate platoons of varying sizes from $\displaystyle 1$ to $\displaystyle 5$ vehicles. The speed limit of road is $\displaystyle 18~m/s$. The maximum speed is set as $\displaystyle 18~m/s$ for platoons going straight,   $\displaystyle 9~m/s$ for left turning and  $\displaystyle 7~m/s$ for right turning platoons. We set the maximum acceleration limit to be $3~m/s^2$ and the minimum deceleration to be $-3~ m/s^2$. We implement the upper-level and the low-level optimal framework in  MATLAB. In the upper-level, we collect the attributes of platoons including link number, lane number, path, current position, current speed, and  number of following vehicles. Then, we compute the time of entry for each platoon into the merging zone. In the low-level, we derive the optimal control input and computes the speed of each platoon based on the time of entry provided by the upper-level framework. Then, the speed of the platoons is updated using COM interface in VISSIM  traffic simulator in real-time.  
 \par   The vehicles form platoons of varying sizes from $1$ to $5$ and are allowed to cross the intersection based on proposed optimal framework (denoted as OC\_Platoon). To validate the effectiveness of the proposed optimal framework, we performed a comparative study with four cases.\\
 \textbf{Case 1:} The platoons of varying sizes from $1$ to $5$ are allowed to cross the intersection based on  FCFS scheduling algorithm and this case is denoted as FCFS\_Platoon. \\
  \textbf{Case 2:}  The vehicles are not allowed to cross the intersection as platoons but they are allowed to cross the intersection one after the other as individual vehicles based on FCFS scheduling algorithm and this case is denoted as FCFS\_ind.\\
 \textbf{Case 3:} The vehicles are not allowed to cross the intersection as platoons but they are allowed to cross the intersection one after the other  based on proposed optimal control algorithm and this case is denoted as OC\_ind.\\  
 \textbf{Case 4:}  The vehicles are given right of way at the intersection based on LQF-MWM scheduling algorithm \cite{maguluri2014stability}. Each lane approaching the intersection is assigned a weight based on their queue length.  We group  lanes of non-conflicting traffic movements  as compatible lanes, i.e., vehicles in compatible lane group  can be given right of way concurrently inside the intersection. For instance (Fig. 1), Lanes 2, 3, 5 \& 6 ,  Lanes 1 \& 4,  Lanes 8, 9, 10, \& 11,  and Lanes 7 \& 12 are compatible lanes.  The weights associated with each lane group is sum of queue lengths of each lane in the group. The weights are calculated for a time interval. In each time interval, the lane group with maximum weight is given right of way inside the intersection.

\begin{figure}[htbp!]
	\centering
	\includegraphics[width=2.8in]{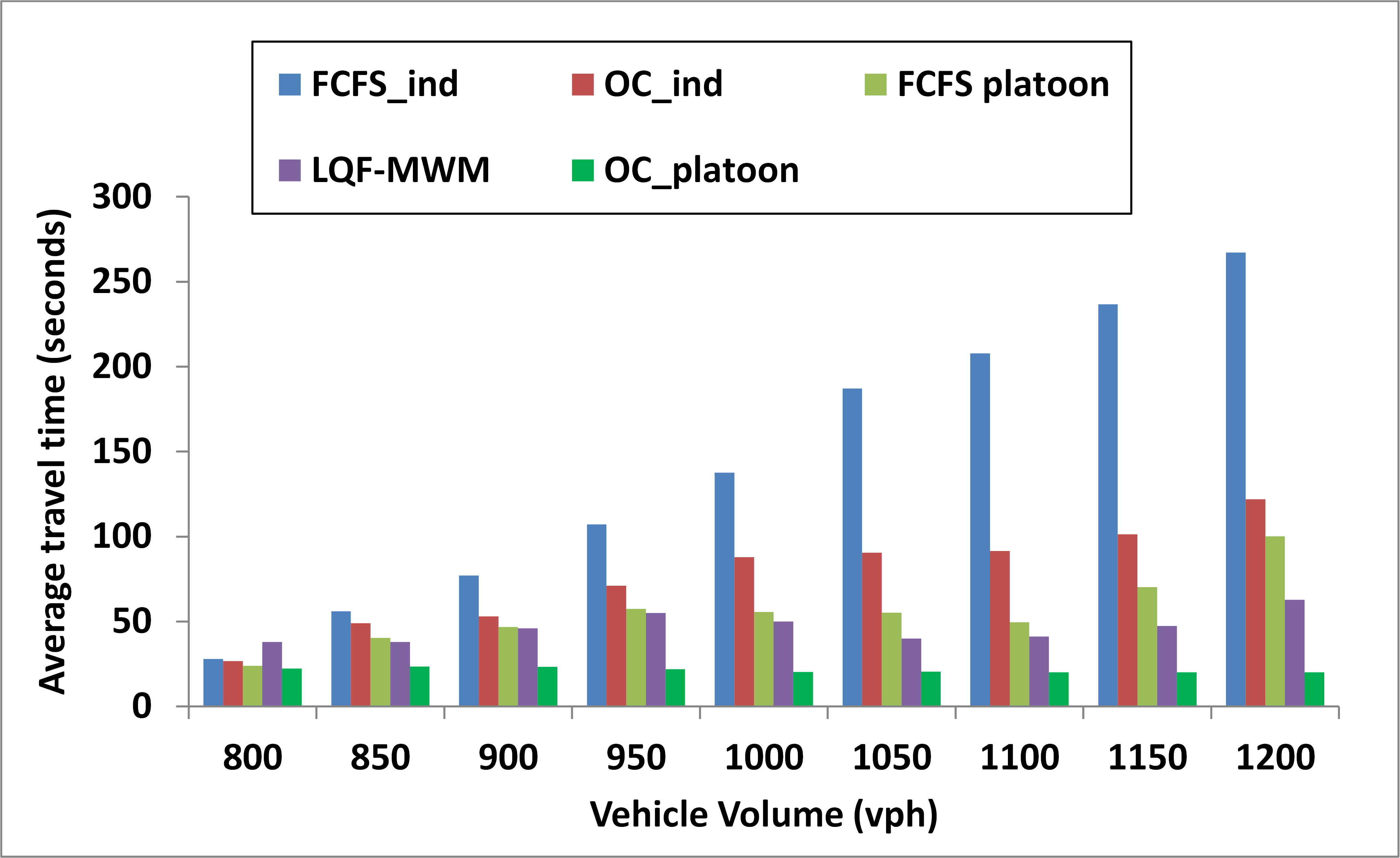}
	\caption{Average travel time.}
	\label{avgTT}	
\end{figure}

\begin{figure}[htbp!]
	\centering
	\includegraphics[width=2.8in]{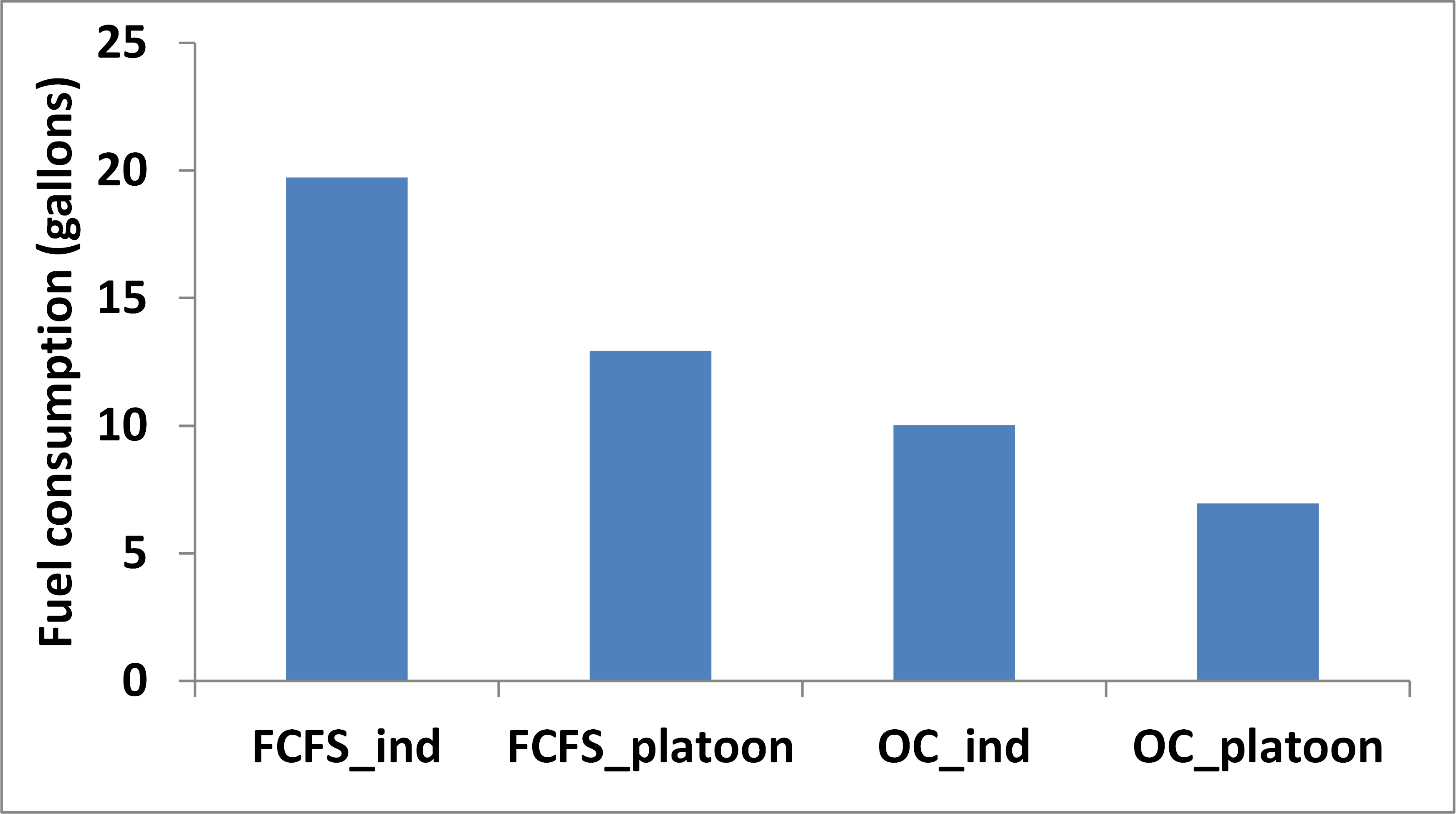}
	\caption{Fuel consumption for each case.}
	\label{fuel1}
\end{figure}

 We ran the simulation for $900$ seconds and collected the evaluation data to compare the performance of the proposed framework in terms of  average travel time  and fuel consumption with FCFS\_Platoon, LQF-MWM, FCFS\_ind and OC\_ind cases.  The  average travel time  and fuel consumption for all the cases  are shown in Figs. \ref{avgTT} and \ref{fuel1}.

 \begin{figure}[htbp!]
 	\centering
 \includegraphics[width=2.8in]{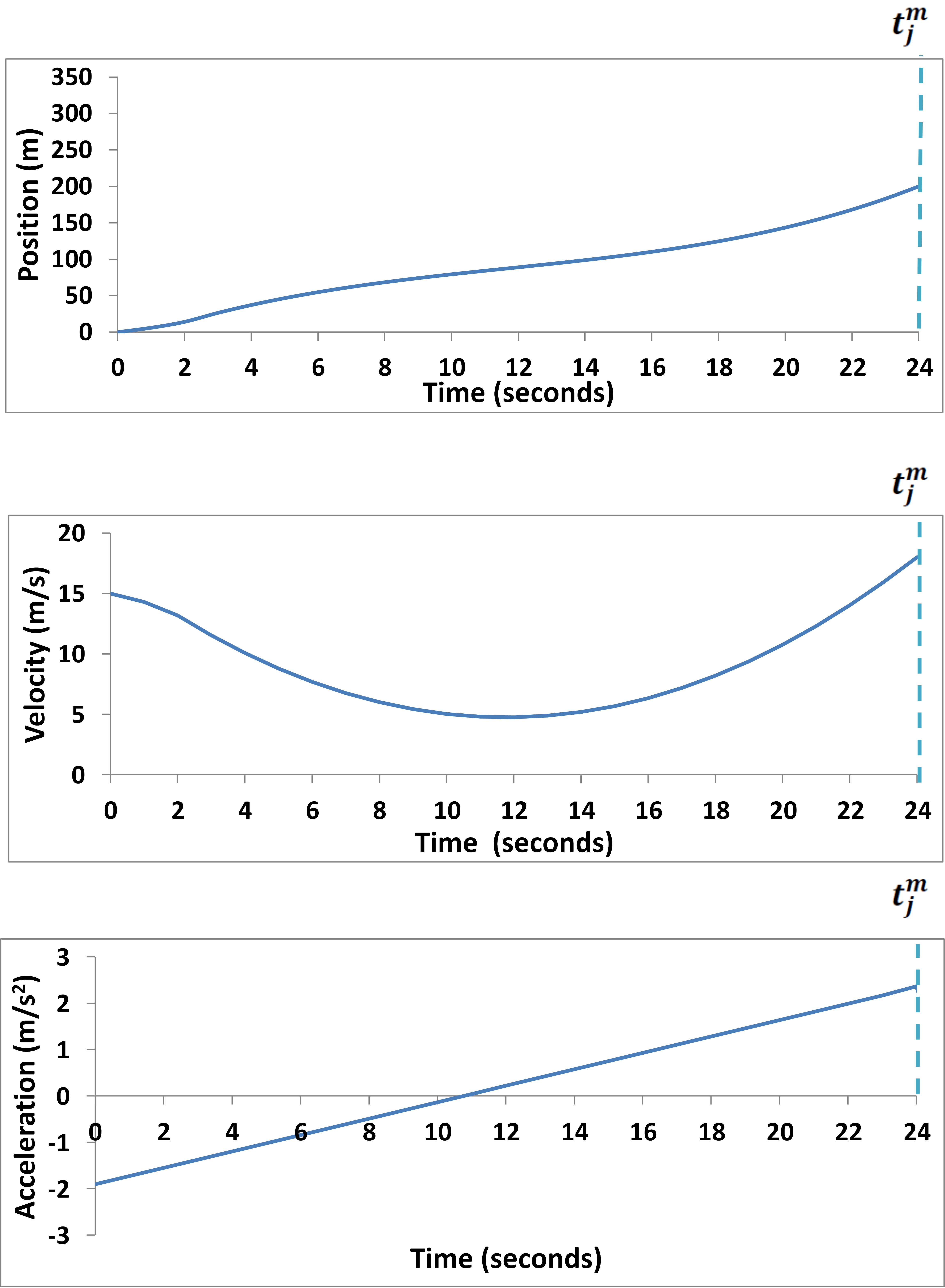}
 	\caption{Position, velocity, and acceleration profile of a platoon.}
 	\label{trajectory}
 \end{figure}
 \par The position, velocity, and acceleration profiles of a platoon from the time of entry inside the scheduling zone to the time at which it exits the merging zone is shown in Fig. \ref{trajectory}. The position profile of the platoon indicates that the platoons enter the merging zone at the time of entry provided by the upper-level optimization framework. The speed profile of the platoon shows that the platoon updates the speed based on their optimal control input (acceleration/deceleration) and exits the merging zone without stopping at the intersection. In case of increasing traffic flow, few platoons stops and wait for the previously entered platoons in the schedule zone to exit the intersection. Then, the platoons find the optimal speed to exit the intersection.

\par The comparison of performance measures of vehicles under OC\_Platoon, FCFS\_Platoon, LQF-MWM, OC\_Ind cases with FCFS\_ind are summarized
in Table \ref{OC}.
\begin{table}
\centering
\caption{Comparison of performance measures of Vehicles under  different cases with FCFS\_Ind case}
\label{OC}
\resizebox{0.7 \columnwidth}{!}{%
\begin{tabular}{cccc}
\hline
Cases & \begin{tabular}[c]{@{}c@{}}Average\\  Travel Time\\   (sec)\end{tabular} & \begin{tabular}[c]{@{}c@{}}Fuel \\ Consumption\\ (gallons)\end{tabular} \\ \hline
OC\_Ind & -46.89\%  & -49.1\% \\
FCFS\_Platoon & -61.74\%  & -52.5\% \\
\multicolumn{1}{l}{LQF-MWM} & -68.71\%  & -22.79\% \\
OC\_Platoon & -84.96\%  & -64.76\% \\ \hline
\end{tabular}%
}
\end{table}
 The OC\_ind case reduced the average travel time by 46.89\% and fuel consumption by 49.1\% over FCFS\_Ind case.  The FCFS\_platoon case reduced the average travel time by 61.74\%  and fuel consumption by 52.5\% over FCFS\_Ind case. The LQF-MWM algorithm reduced the the average travel time by 68.71\% and fuel consumption by 22.79\% over FCFS\_Ind case. The LQF-MWM algorithm resulted in high fuel consumption since vehicles are completely stopped at stop line and given right of way based on the queue length of compatible lane groups. The OC\_platoon case reduced the average travel time by 84.96\% and fuel consumption by 64.76\% over FCFS\_Ind case. 
 
 \par
\begin{figure}[htbp!]
	\centering
	\includegraphics[width=2.7in]{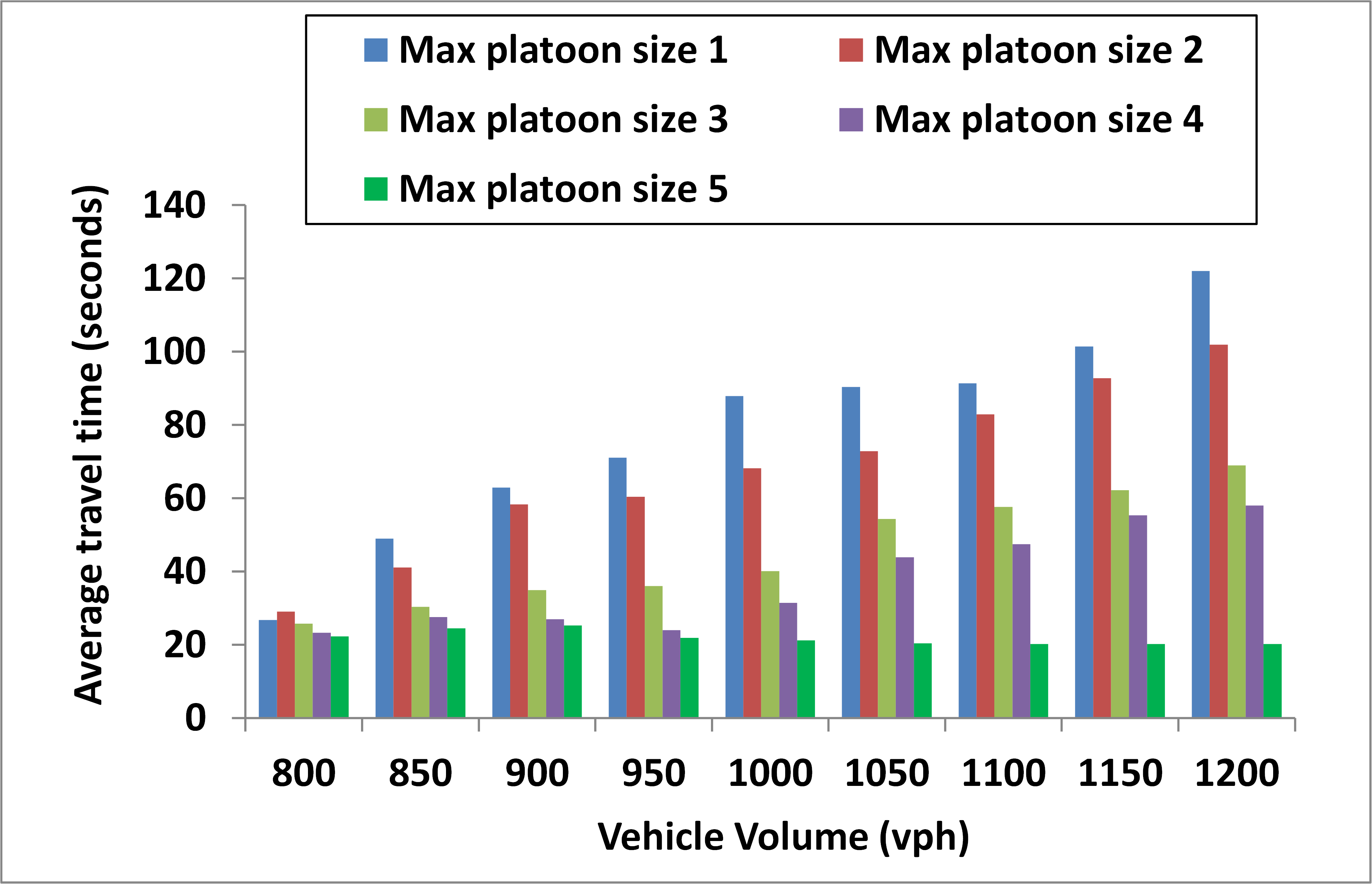}
	\caption{Average travel time of vehicles.}
	\label{avgTTP}
\end{figure}

\begin{figure}[htbp!]
	\centering
	\includegraphics[width=2.7in]{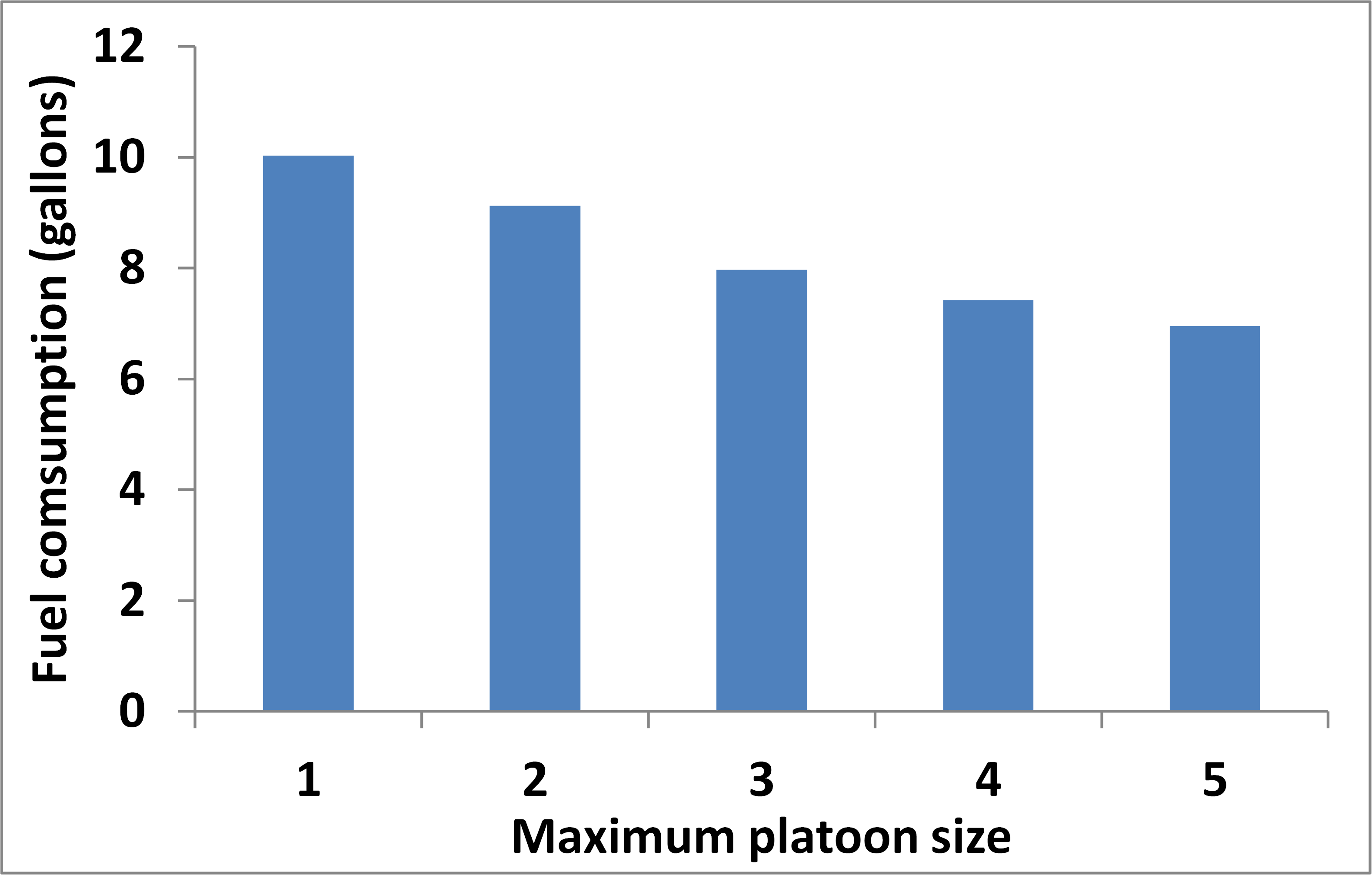}
	\caption{Fuel consumption of vehicles in the network.}
	\label{fuel2}
\end{figure}

 We perform a simulation study to evaluate the effect of platooning while crossing the intersection. We consider five cases where the vehicles are allowed to form platoons of maximum size 1, 2, 3, 4, and 5, respectively and allowed to cross the intersection based on proposed optimal control algorithm. The average travel time and fuel consumption for all the cases of varying platoon sizes  are shown in Figs. \ref{avgTTP} and \ref{fuel2}. The comparison of performance measures of vehicles under different maximum platoon sizes cases with maximum platoon size 1  are summarized
in Table \ref{diffsize}. 
\begin{table}
\centering
\caption{Comparison of performance measures of different maximum platoon sizes with maximum platoon size 1}
\label{diffsize}
\resizebox{0.6\columnwidth}{!}{%
\begin{tabular}{cccc}
\hline
\begin{tabular}[c]{@{}c@{}}Platoon\\ Size\end{tabular} & \begin{tabular}[c]{@{}c@{}}Average\\  Travel Time\\   (sec)\end{tabular}  & \begin{tabular}[c]{@{}c@{}}Fuel \\ Consumption\\ (gallons)\end{tabular} \\ \hline
2 & -13.56\%  & -9.05\% \\
3 & -41.61\% & -20.6\% \\
4 & -51.9\%  & -26.03\% \\
5 & -71.68\%  & -30.72\% \\ \hline
\end{tabular}%
}
\end{table}
The scenario in which the vehicles are allowed to form platoon of  maximum size  2 reduced the average travel time by 13.56\%  and fuel consumption by 9.05\% when compared with scenario in which the maximum platoon size is 1. The scenario in which the vehicles are allowed to form platoon of  maximum size 3 reduced the average travel time by 41.61\% and fuel consumption by 20.6\% when compared with scenario in which platoon size is 1.  The scenario in which the vehicles are allowed to form platoon of  maximum size 4 reduced the average travel time by 51.9\% and fuel consumption by 26.03\% when compared with scenario in which platoon size is 1. The scenario in which the vehicles are allowed to form platoon of  maximum size 5 reduced the average travel time by 71.68\% and fuel consumption by 30.72\% when compared with scenario in which platoon size is 1.

\section{Concluding Remarks and Future Work}  
In this paper, we investigated the problem of optimal coordination of platoons of CAVs at a signal-free intersection. We developed a decentralized, two-level optimal framework for coordinating the platoons with the objective to minimize travel delay and fuel consumption. In the upper level, we presented an optimization framework to reduce the delay of the platoons at an intersection. In the low level, we presented  a time and energy optimal control problem, and derived analytical solutions that provided optimal control input to the platoons. We performed a comparative study of the proposed framework with the FCFS and LQF-MWM scheduling algorithms. The simulation analysis showed that the proposed framework significantly reduces the travel time and fuel consumption of the platoons at the intersection. We also investigated the effect of platooning by considering platoons of varying sizes while crossing the intersection.

 Ongoing efforts consider lane changes of platoons at an intersection. The proposed approach assumes that vehicles form platoons before entering the schedule zone and restricts the vehicle from one platoon to join other platoon inside the schedule zone. Future research should  consider the formation of platoons inside the schedule zone along with stability of the platoons, and  extend the proposed framework for a mixed environment of human-driven vehicles and CAVs at different penetration rates. 
\vspace{-4mm}

\section*{Acknowledgment}
The authors would like to thank Behdad Chalaki, A M Ishtiaque Mahbub and L. E. Beaver for the technical discussions.
\vspace{-4mm}
\bibliographystyle{IEEEtran}        % Include this if you use bibtex 
\bibliography{IEEEabrv,ref.bib}

%\vskip 0pt plus -1fil
\vspace{-14mm}
\begin{IEEEbiography}[{\includegraphics[width=1.1in,height=2in,clip,keepaspectratio]{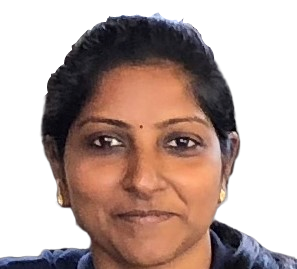}}]
	{Sharmila Devi Kumaravel}
 received the B.E degree in Electronics and Instrumentation Engineering from Bharathidasan University, Trichy in 2004 and M.E degree in Power Electronics and Drives from Anna University, Chennai in 2008. She is currently working for her Ph.D. degree with the Department of Instrumentation and Control Engineering, National Institute of Technology, Tiruchirappalli, India. She was a  visiting research scholar in the Department of Mechanical Engineering at the University of Delaware. Her research interests include  Intelligent transportation systems, Control Engineering and Graph theory.
\end{IEEEbiography}
\vspace{-10mm}
%\vskip 0pt plus -1fil
\begin{IEEEbiography}[{\includegraphics[width=1.1 in,height=2.2in,clip,keepaspectratio]{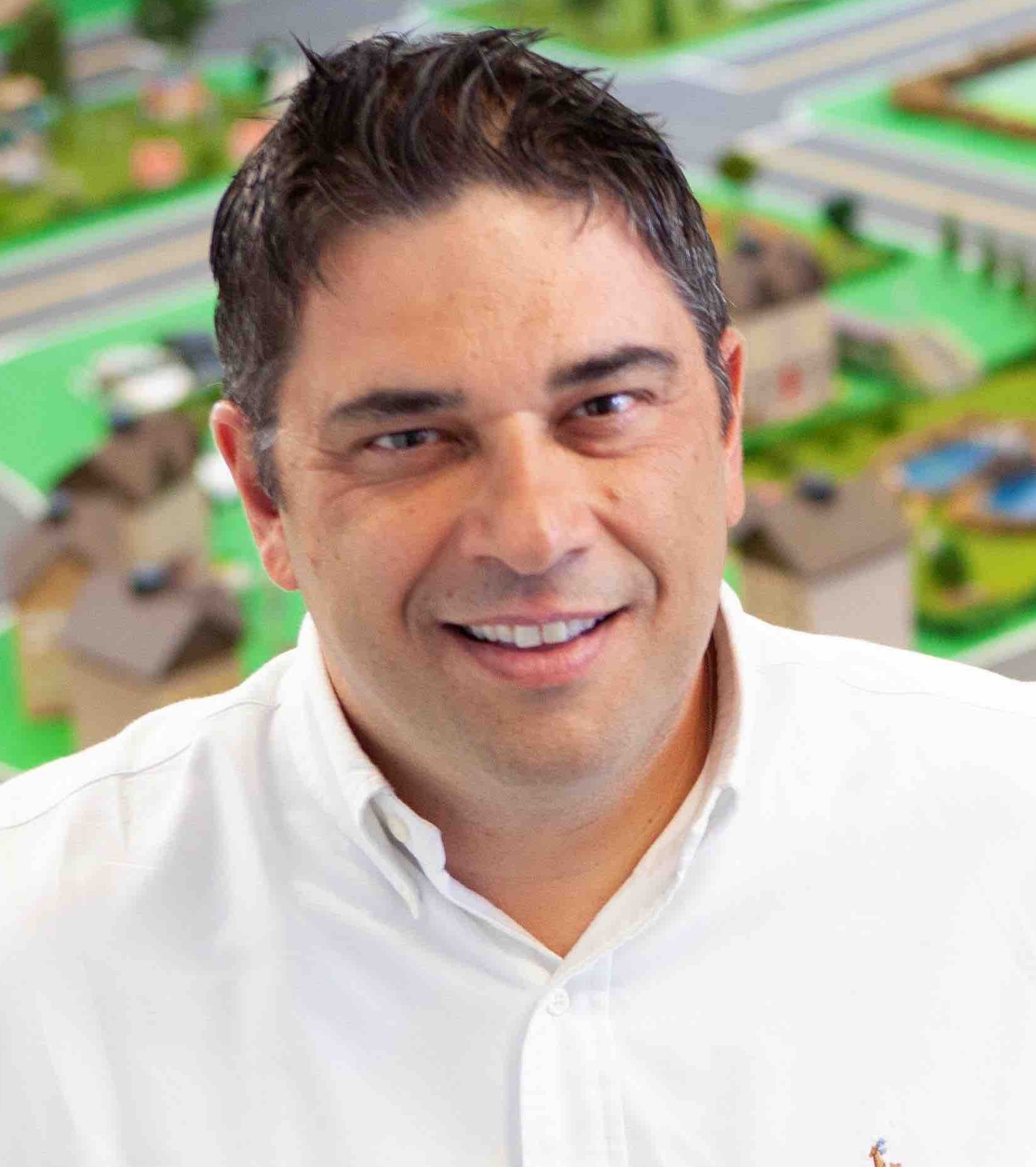}}]
	{Andreas A. Malikopoulos}
   (M2006, SM2017) received the Diploma in mechanical engineering from the National Technical University of Athens, Greece, in 2000. He received M.S. and Ph.D. degrees from the department of mechanical engineering at the University of Michigan, Ann Arbor, Michigan, USA, in 2004 and 2008, respectively. 
	He is the Terri Connor Kelly and John Kelly Career Development Associate Professor in the Department of Mechanical Engineering at the University of Delaware (UD), the Director of the Information and Decision Science (IDS) Laboratory, and the Director of the Sociotechnical Systems Center. Before he joined UD, he was the Deputy Director and the Lead of the Sustainable Mobility Theme of the Urban Dynamics Institute at Oak Ridge National Laboratory, and a Senior Researcher with General Motors Global Research \& Development. His research spans several fields, including analysis, optimization, and control of cyber-physical systems; decentralized systems; stochastic scheduling and resource allocation problems; and learning in complex systems. The emphasis is on applications related to smart cities, emerging mobility systems, and sociotechnical systems. He has been an Associate Editor of the IEEE Transactions on Intelligent Vehicles and IEEE Transactions on Intelligent Transportation Systems from 2017 through 2020. He is currently an Associate Editor of Automatica and IEEE Transactions on Automatic Control. He is a member of SIAM, AAAS, and a Fellow of the ASME.
\end{IEEEbiography}
\vspace{-10mm}
%\vskip 0pt plus -1fil
\begin{IEEEbiography}[{\includegraphics[width=1.1in,height=1.5in,clip,keepaspectratio]{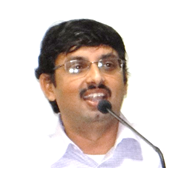}}]{Ramakalyan Ayyagari}
 is a professor in Instrumentation \& Control Engineering Dept., National Institute of Technology, Tiruchirappalli, India. He is deeply interested in looking into computational problems that arise out of the algebra and graphs in control theory and applications. Of particular interest are the NP-hard problems and the Randomized Algorithms. He is a senior member of IEEE and member of SIAM. He was the founding secretary and past President of Automatic Control and Dynamic Optimization Society (ACDOS), the Indian NMO of the IFAC, through which he passionately contributes to controls education in India.	
\end{IEEEbiography}

 \end{document}